\documentclass[preprint,11pt]{elsarticle}

\usepackage{lmodern}
\usepackage{amsmath,amssymb,amsthm,mathtools}
\usepackage{booktabs}
\usepackage{microtype}
\usepackage[hidelinks]{hyperref}
\hypersetup{
  pdftitle={The Ternary and Quaternary Cases of the \{2,4\}-Distance Code Conjecture},
  pdfauthor={Valjakas Djaljapayan; Wei-Hsuan Yu; Yen-Chi Roger Lin}
}
\usepackage{enumitem}
\setlist[enumerate]{leftmargin=*,topsep=4pt,itemsep=2pt}

\newtheorem{theorem}{Theorem}[section]
\newtheorem{lemma}[theorem]{Lemma}

\newtheorem{corollary}[theorem]{Corollary}
\theoremstyle{remark}

\newcommand{\supp}{\operatorname{supp}}
\newcommand{\wt}{\operatorname{wt}}

\journal{Journal of Combinatorial Theory, Series A}

\begin{document}
\begin{frontmatter}

\title{The Ternary and Quaternary Cases of the $\{2,4\}$-Distance Code Conjecture}
\author[ncu]{Valjakas Djaljapayan}
\ead{115281004@cc.ncu.edu.tw}

\author[ncu]{Wei-Hsuan Yu}
\ead{whyu@math.ncu.edu.tw}

\author[ntnu]{Yen-Chi Roger Lin\corref{cor1}}
\ead{yclin@math.ntnu.edu.tw}

\cortext[cor1]{Corresponding author.}

\affiliation[ncu]{
  organization={Department of Mathematics, National Central University},
  addressline={No. 300, Zhongda Road, Zhongli District},
  city={Taoyuan City},
  postcode={320317},
  country={Taiwan}
}

\affiliation[ntnu]{
  organization={Department of Mathematics, National Taiwan Normal University},
  addressline={No. 88, Section 4, Ting Chou Road},
  city={Taipei},
  postcode={116},
  country={Taiwan}
}

\begin{abstract}
For an integer $q\ge2$, let $A_q(n,\{2,4\})$ be the largest cardinality of
an unrestricted $q$-ary code of length $n$ in which any two distinct words
have Hamming distance $2$ or $4$.  We prove that
\[
 A_q(n,\{2,4\})=1+\binom n2
 \qquad(q\in\{3,4\},\ n\ge6),
\]
thereby settling the two remaining alphabet sizes in the conjecture of
Boyvalenkov, Delchev, Zinoviev, and Zinoviev; the binary case was settled
previously.  The quaternary result implies the ternary result by alphabet
inclusion.  For $n\ge16$, the quaternary upper bound is proved
combinatorially.  After normalizing one word to zero, the weight-four
supports form a $2$-intersecting family; exact support intersections force
symbol synchronization and restrict support multiplicities.  We combine
the nontrivial intersection theorem with matching arguments and a
three-coordinate projection lemma.  The finite lengths $6\le n\le13$
are handled by a complete Hamming-isometry reduction to exact clique
decisions, and $n=14,15$ by a smaller exact weighted-clique relaxation.
All finite decisions are reproduced by two standalone integer-only C++17
verifiers provided as supplementary source files.
\end{abstract}

\end{frontmatter}

\section{Introduction}

For an alphabet $\Sigma$ with $|\Sigma|=q$, Hamming distance on $\Sigma^n$ counts the
coordinates on which two words differ.  If $D\subseteq\{1,\ldots,n\}$,
write $A_q(n,D)$ for the largest size of a subset $C\subseteq\Sigma^n$ for
which $d(x,y)\in D$ whenever $x\ne y$.  In particular, codes need not
be linear, and the alphabet need not be a field.

The classical coding-theory background is developed in
\cite{Hamming50,MS77,vanLint99,HP03}.
The study of prescribed-distance codes is related to Delsarte's
association-scheme inequalities \cite{Delsarte73,BannaiIto84,BCN89},
to semidefinite refinements \cite{Sch05,GST06,BGSV12}, and to
intersection theorems for uniform set systems
\cite{EKR61,HM67,Frankl78,FW81,Wilson84,AK96,AK97,FT18,GodsilMeagher16}.
Earlier work on codes with few distances includes
\cite{EFIN87,CK86,BDZZ21,LRV23}.  The analogous Euclidean and spherical
problems provide related techniques and bounds
\cite{DGS77,LRS77,BBS83,Musin09,GY18}.

The construction
\begin{equation}\label{eq:standard}
 C_{\mathrm{std}}=\{0\}\cup\{w_{ij}:1\le i<j\le n\},
\end{equation}
where $w_{ij}$ has a distinguished nonzero symbol in positions $i,j$ and
zero elsewhere, has pairwise distances in $\{2,4\}$.  Thus
\begin{equation}\label{eq:lower}
 A_q(n,\{2,4\})\ge1+\binom n2
 \qquad(q\ge2).
\end{equation}
Boyvalenkov, Delchev, Zinoviev, and Zinoviev
\cite[Conjecture~1(i)]{BDZZ21} conjectured that equality holds for
$q\in\{2,3,4\}$ and every $n\ge6$.  The binary case follows from the
work of Barg, Glazyrin, Kao, Lai, Tseng, and Yu \cite{Barg24}.
We establish the two remaining cases.

\begin{theorem}\label{thm:main}
For each integer $n\ge6$ and each $q\in\{3,4\}$,
\[
 A_q(n,\{2,4\})=1+\binom n2.
\]
\end{theorem}

It suffices to prove the assertion for $q=4$.  Indeed, if $q\le r$,
every $q$-ary code embeds isometrically in an $r$-ary Hamming space, so
\begin{equation}\label{eq:alphabet-monotone}
 A_q(n,D)\le A_r(n,D).
\end{equation}
The construction \eqref{eq:standard} then gives the same lower bound for
$q=3$.  This observation avoids repeating the quaternary argument in the
ternary setting.

Our proof uses a support--symbol interaction specific to the forbidden
odd distances.  Once a word is normalized to zero, every other word
has weight two or four.  The family of weight-four supports is
$2$-intersecting; if two of those supports meet in precisely two
coordinates, the symbols on their intersection must agree.  Consequently,
the usual unweighted intersection bounds can be strengthened by local
multiplicity constraints.  If all weight-four supports contain a common
two-set, the remaining coordinate pairs form a graph whose matching
number governs the argument.  The only additional ingredient required
for the four-symbol alphabet is a bound on projections to three
coordinates.

These structural arguments prove the claim for $n\ge16$.  For the ten
remaining lengths we use exact, reproducible finite decisions, described
mathematically in Section~\ref{sec:finite}.  Such computer-assisted
methods are standard in extremal clique calculations
\cite{CarraghanPardalos90,TomitaKameda07}; the present verification uses
neither floating-point optimization nor an unverified heuristic.
The finite computations and their supporting code are logically separate
from the infinite-range proof.

\section{Normalization and the weight-two layer}

We now work in an alphabet $\Sigma=\{0,1,2,3\}$.  Arbitrary permutations of
coordinates and independent permutations of the four symbols in each
coordinate preserve Hamming distance.  We may therefore assume that
$0\in C$ and write
\begin{equation}\label{eq:layers}
 C=\{0\}\sqcup C_2\sqcup C_4,
\end{equation}
where $C_i$ consists of the words of weight $i$.
For a word $x$, its \emph{support} is
$\supp(x)=\{j:x_j\ne0\}$, and $\wt(x)=|\supp(x)|$.
For a support $S$ occurring in a weight layer, its \emph{multiplicity}
is the number of words in that layer with support $S$.

\begin{lemma}[One-coordinate synchronization]\label{lem:one}
Let $x\in C_2$ and $y\in C_2\cup C_4$.  If
$\supp(x)\cap\supp(y)=\{i\}$, then $x_i=y_i$.
\end{lemma}
\begin{proof}
Let $\delta=1$ if $x_i\ne y_i$, and $\delta=0$ otherwise.
Since $i$ is the only common support coordinate,
\[
 d(x,y)=\wt(x)+\wt(y)-2+\delta.
\]
This is $2+\delta$ if $y$ has weight two, and $4+\delta$ if $y$ has
weight four.  The distances $3$ and $5$ are forbidden, so $\delta=0$.
\end{proof}

For a code $B\subseteq\Sigma^N$ consisting of weight-two words with
mutual distances in $\{2,4\}$, let $\mathcal E$ be its family of distinct
supports.  A support is \emph{isolated} if it is disjoint from every
other member of $\mathcal E$.  Distinct words on the same support must
differ in both coordinates, so each support has multiplicity at most
three.  If a support meets another member of $\mathcal E$,
Lemma~\ref{lem:one} fixes the symbol at their common coordinate;
its multiplicity is therefore one.  If $r$ supports are isolated,
they occupy $2r$ coordinates, and every other support avoids these
coordinates.  Consequently,
\begin{equation}\label{eq:pure-w2}
 \begin{aligned}
 |B|&\le|\mathcal E|+2r
      \le3r+\binom{N-2r}{2}\\
    &=\binom N2-r(2N-2r-4)
      \le\binom N2 \qquad(N\ge3),
 \end{aligned}
\end{equation}
where the last inequality follows from $r\le\lfloor N/2\rfloor$ and
$2N-2r-4\ge2\lceil N/2\rceil-4\ge0$.

\begin{lemma}[Two useful bounds for $C_2$]\label{lem:C2-bounds}
Assume $C_4\ne\varnothing$.
Then
\begin{equation}\label{eq:C2-general}
 |C_2|\le4n-2.
\end{equation}
If $C_4$ contains two distinct supports whose intersection has size two,
then
\begin{equation}\label{eq:C2-exact}
 |C_2|\le2n+9.
\end{equation}
\end{lemma}
\begin{proof}
Fix a word $z\in C_4$ and put $S=\supp(z)$.  Every weight-two
support meets $S$, for otherwise its distance from $z$ would be six.
There are
\[
 \binom n2-\binom{n-4}{2}=4n-10
\]
such two-supports.  Their isolated members are pairwise disjoint,
so there are at most four of them.  By the multiplicity observation above,
each isolated support contributes at most two words beyond a baseline of
one word per support.  This proves \eqref{eq:C2-general}.

For the second bound, let the two four-supports be
$S=P\sqcup A$ and $T=P\sqcup B$, where $P,A,B$ are disjoint two-sets.
Every weight-two support meets both $S$ and $T$; inclusion--exclusion
gives at most
\[
 \binom n2-2\binom{n-4}{2}+\binom{n-6}{2}=2n+1
\]
possible supports.  An isolated one meeting $P$ uses at least one of its
two coordinates, so at most two such supports exist.  An isolated
support avoiding $P$ is a pair with one endpoint in $A$ and one in $B$;
pairwise disjoint pairs of this type form a matching in $K_{2,2}$, and
there are at most two.  Thus there are at most four isolated supports.
Adding at most eight extra words proves \eqref{eq:C2-exact}.
\end{proof}

\section{Local multiplicities on four-supports}

Let $\mathcal H$ denote the set of distinct supports of words in $C_4$.
For $S\in\mathcal H$, write
\[
 \mu(S)=|\{x\in C_4:\supp(x)=S\}|.
\]
Only these two notations, $\mathcal H$ and $\mu(S)$, will be needed
for the weighted support family.

\begin{lemma}[A single four-support]\label{lem:local9}
Every $S\in\mathcal H$ satisfies $\mu(S)\le9$.
\end{lemma}
\begin{proof}
Restrict the words on $S$ to their four nonzero coordinates, regarding
the three nonzero symbols as an alphabet $\Omega$ of size three.  Let
$M=\mu(S)$ and put
\[
 a_j=\frac1M|\{(x,y):x,y\text{ are in the restricted code},
 \ d(x,y)=j\}|.
\]
Thus $a_0=1$, $a_1=a_3=0$, and $M=1+a_2+a_4$.
For $0\le j\le4$ define the ternary Krawtchouk coefficient
\[
 K_3(j)=\sum_{r=0}^3(-1)^r2^{3-r}
 \binom jr\binom{4-j}{3-r}.
\]
The standard positivity inequality gives
$\sum_{j=0}^4a_jK_3(j)\ge0$.
For completeness, choose a real orthonormal basis
$1,\phi_1,\phi_2$ of functions on $\Omega$ under the uniform inner
product.  Completeness yields
\[
 \sum_{i=1}^2\phi_i(u)\phi_i(v)
 =3\mathbf1_{\{u=v\}}-1.
\]
For each three-subset $I$ of the four coordinates, form all products
$\prod_{i\in I}\phi_{b_i}(x_i)$ with $b_i\in\{1,2\}$.
Expanding the sum of the squares of their totals over the restricted
code, and dividing by $M$, gives exactly
$\sum_j a_jK_3(j)$.  This proves the required positivity without a
linearity assumption.

Direct evaluation gives $K_3(0)=32$ and
$K_3(2)=K_3(4)=-4$.  Hence
$0\le32-4a_2-4a_4$, so $M=1+a_2+a_4\le9$.
\end{proof}

\begin{lemma}[Exact two-overlaps]\label{lem:overlaps}
Let $S,T\in\mathcal H$ be distinct and satisfy $|S\cap T|=2$.
Then $\mu(S),\mu(T)\le3$.  Moreover, if another support
$U\in\mathcal H$ satisfies $|S\cap U|=2$ and
$S\cap U\ne S\cap T$, then $\mu(S)=1$.
\end{lemma}
\begin{proof}
Put $P=S\cap T$.  If $x$ is supported on $S$ and $y$ on $T$, their
four coordinates in the symmetric difference $S\triangle T$
already contribute distance four.  Thus $x,y$ must agree on $P$.
In particular, all words on $S$ have the same two symbols on $P$.
Two distinct such words can differ only on $S\setminus P$, and they
must differ in both remaining coordinates.  Projection to either one
is injective into the three nonzero symbols, giving $\mu(S)\le3$.
Interchanging $S,T$ gives the other bound.

If $U$ gives a distinct exact-overlap pair, then comparisons with words
on $T$ and $U$ fix the symbols of every word on $S$ on the union of two
distinct two-subsets of $S$.  This union has size at least three.
Two different words on $S$ would then be at distance at most one, an
impossibility.  Hence $\mu(S)=1$.
\end{proof}

\begin{lemma}[The support family is $2$-intersecting]
\label{lem:two-intersecting}
For all distinct $S,T\in\mathcal H$ one has $|S\cap T|\ge2$.
\end{lemma}
\begin{proof}
Choose words $x,y$ on $S,T$ and put $t=|S\cap T|$.  With $\delta$
the number of unequal nonzero symbols on the common support,
\[
 d(x,y)=8-2t+\delta\le4.
\]
Since $\delta\ge0$, we obtain $t\ge2$.
\end{proof}

\section{Three-coordinate projections}

The following two elementary lemmas isolate the only genuinely
quaternary local calculation.

\begin{lemma}[Distance-two sets in three coordinates]
\label{lem:distance2-three}
If $B\subseteq\Omega^3$, where $|\Omega|=3$, has
$d(u,v)=2$ for every pair of distinct words, then $|B|\le4$.
\end{lemma}
\begin{proof}
If $|B|\le1$ there is nothing to prove.  Otherwise normalize two
words to $000$ and $110$ by coordinatewise symbol permutations.
A third word at distance two from both has one of the forms
\[
 220,\qquad 10c,\qquad 01c\quad(c\in\{1,2\}).
\]
At most one word of each of the last two forms can occur.  The word
$220$ cannot occur together with either form, because its distance
from either is three.  Consequently $|B|\le4$.
\end{proof}

\begin{lemma}[Three-coordinate projection]\label{lem:projection}
Suppose $B\subseteq\Omega^4$, $|\Omega|=3$, has all mutual distances
in $\{2,4\}$.  Project onto any three coordinates.  If no two
projected words have distance three, then $|B|\le4$.
\end{lemma}
\begin{proof}
The projection is injective: otherwise two distinct original words
would differ only in the omitted coordinate.  Denote the projected set
by $Y$, and for $u\in Y$ let $f(u)$ be the omitted symbol.
For distinct $u,v\in Y$ the projected distance is one or two, and
\begin{equation}\label{eq:projection-rules}
 d(u,v)=1\ \Longrightarrow\ f(u)\ne f(v),\qquad
 d(u,v)=2\ \Longrightarrow\ f(u)=f(v).
\end{equation}
Make a graph on $Y$ by joining words at projected distance two.
The value of $f$ is constant on each connected component.  Consequently,
all distinct vertices within a component have projected distance two,
and Lemma~\ref{lem:distance2-three} bounds each component by four.

If there is just one component, this proves the result.  Otherwise
words from different components have projected distance one.
If every component is a singleton, all words in $Y$ have pairwise
distance one; fixing two shows that all words differ only in the same
coordinate, so $|Y|\le3$.

Finally, suppose one component contains two words; normalize them to
$000,110$.  A word outside this component must have distance one from
both, and hence is $100$ or $010$.  If exactly one of these occurs,
say $100$, any additional word in the component must be of the form
$10c$ with $c\ne0$; at most one such word can occur because two of
them would be at distance one.  If both $100$ and $010$ occur, there
is no additional word in the component.  Thus $|Y|\le4$ in all cases.
\end{proof}

\begin{lemma}[The antipodal shell]\label{lem:shell}
Let $u,v\in\Omega^3$ have distance three, where $|\Omega|=3$.
Set
\[
 X(u,v)=\{w\in\Omega^3:d(w,u)=d(w,v)=2\}.
\]
Join two elements of $X(u,v)$ when their distance is two.
The resulting graph is triangle-free.
\end{lemma}
\begin{proof}
Normalize $u=000$ and $v=111$.  The six elements of $X(u,v)$ are the
permutations of $012$.  Each such element has exactly three neighbors
at distance two, and those three neighbors are pairwise at distance
three.  For example, the neighbors of $012$ are $021,102,210$.
Thus no vertex belongs to a triangle.
\end{proof}

\begin{lemma}[Cross-families in a triangle-free graph]
\label{lem:cross}
Let $G$ be triangle-free, and let nonempty
$B_1,\ldots,B_r\subseteq V(G)$ satisfy the following conditions:
for $i\ne j$, every $x\in B_i$, $y\in B_j$ are equal or adjacent;
$|B_i|\le6$; and each independent $B_i$ has size at most three.
Then
\begin{equation}\label{eq:cross}
 \sum_{i=1}^r|B_i|\le2r+4.
\end{equation}
\end{lemma}
\begin{proof}
If some $B_i$ contains an edge $ab$, every vertex of any other
$B_j$ must be $a$ or $b$; a different vertex adjacent to both would
create a triangle.  Hence the total is at most $6+2(r-1)$.

Otherwise every $B_i$ is independent.  If at most one of these sets
has size three, the sum is at most $3+2(r-1)\le2r+4$.
If two have size three, call them $B_1,B_2$.
Their intersection has size at most one, since two distinct common
vertices would be nonadjacent within $B_1$ but required to be
adjacent by cross-compatibility.  Choose distinct
$a\in B_1$ and $b\in B_2$, which are adjacent.  Every vertex in a
third set belongs to $\{a,b\}$ by triangle-freeness, and that third
set, being independent, contains at most one of them.  Thus the total
is at most $3+3+(r-2)\le2r+4$ (the case $r=2$ is immediate).
\end{proof}

\begin{lemma}[Common three-core bound]\label{lem:common3}
Suppose $\mathcal T\subseteq\mathcal H$ consists of $s$ distinct
supports all containing the same three-set.  Then
\begin{equation}\label{eq:common3}
 \sum_{S\in\mathcal T}\mu(S)\le\max\{4s,\,2s+11\}.
\end{equation}
\end{lemma}
\begin{proof}
Let the common three-set be $Q$.  Every support in $\mathcal T$ has
the form $Q\cup\{v\}$.  Projection of the words on any one support
to $Q$ is injective, because two distinct words with the same
projection would have distance one.

If every support has multiplicity at most four, the sum is at most
$4s$.  Otherwise choose a support $S_0$ with $\mu(S_0)\ge5$.
For $s=1$, Lemma~\ref{lem:local9} gives
$\mu(S_0)\le9\le13$, as needed.  Assume $s\ge2$.
By Lemma~\ref{lem:projection}, the projected words on $S_0$
contain a pair $u,v$ at distance three.

For a different support $T\in\mathcal T$, let $B_T$ be the set of
three-coordinate projections of its words.  The two distinct leaf
coordinates of $S_0$ and $T$ contribute distance two.  Thus every
$w\in B_T$ has $d(w,u),d(w,v)\in\{0,2\}$; neither can be
zero since $d(u,v)=3$.  Hence
\[
 B_T\subseteq X(u,v).
\]
Label each vertex $w\in B_T$ by the nonzero symbol in the leaf
coordinate of the corresponding word on $T$.  If two projected
vertices are adjacent in the distance-two shell graph, their leaf
symbols agree, since otherwise the original distance would be three.
Hence the leaf symbol is constant on each connected component of the
induced shell graph.  Any two vertices in one component must have
projected distance two, since they have the same leaf symbol and the
full distance is either two or four; the projected distance is at most
three.  Therefore each component is a clique and has size at most two
by Lemma~\ref{lem:shell}.  Different components have different leaf
symbols: two vertices of different components with equal leaf symbols
would have projected distance two and would be adjacent.
There are three nonzero leaf symbols, so $|B_T|\le6$; moreover, if
$B_T$ is independent in the shell graph, then $|B_T|\le3$.

For distinct supports $T,U\ne S_0$, their two leaf coordinates
contribute distance two.  Their projected sets therefore satisfy the
cross-condition of Lemma~\ref{lem:cross}: each pair of projected words
is equal or is at distance two.  Applying that lemma to the
$s-1$ sets $B_T$ gives
\[
 \sum_{T\ne S_0}\mu(T)\le2(s-1)+4=2s+2.
\]
As $\mu(S_0)\le9$, the full sum is at most $2s+11$.
\end{proof}

\begin{corollary}\label{cor:three-small}
If at most three distinct four-supports occur and every pair intersects
in three coordinates, then their total multiplicity is at most $21$.
If all four-supports lie within a fixed five-set, then $|C_4|\le33$.
\end{corollary}
\begin{proof}
If every multiplicity is at most four, the bounds are $12$ and $20$,
respectively.  Otherwise fix a support $S_0$ of multiplicity at least
five.  For each other support $T$, the same projection and antipodal
shell argument in the proof of Lemma~\ref{lem:common3}, using
$S_0\cap T$ as the three-coordinate projection, gives $\mu(T)\le6$.
Thus the totals are at most $9+2\cdot6=21$ and
$9+4\cdot6=33$, respectively.  Distinct four-subsets of a five-set
always intersect in three coordinates.
\end{proof}

\section{Intersection and matching structure}

We use the following specialization of the complete nontrivial
intersection theorem of Ahlswede and Khachatrian \cite{AK96}:
\begin{equation}\label{eq:AK}
 \text{if $n\ge10$ and $\mathcal H$ is $2$-intersecting with
 no common two-set,}\quad |\mathcal H|\le4n-15.
\end{equation}
The cited theorem, with intersection parameter $t=2$ and block size
$k=4$, yields $4\binom{n-4}{1}+\binom{n-4}{0}=4n-15$.
For clarity, a \emph{common two-set} means a fixed two-element set
contained in every member of $\mathcal H$.

We also need two familiar elementary structural facts.

\begin{lemma}[Two intersection and matching facts]\label{lem:structure}
\begin{enumerate}[label=\textup{(\roman*)}]
\item If a family of distinct four-sets has pairwise intersections of
size at least three, then either every set contains a common three-set,
or all sets are four-subsets of a fixed five-set.
\item If a simple graph has no two disjoint edges, then its edges have
a common vertex or the graph is contained in a triangle.
\end{enumerate}
\end{lemma}
\begin{proof}
For (i), choose distinct $S=P\cup\{a\}$ and $T=P\cup\{b\}$ with
$|P|=3$.  If some third set does not contain $P$, its intersections
with $S,T$ force it to be $(P\setminus\{p\})\cup\{a,b\}$ for some
$p\in P$.  Any further member, compared with these three sets, must
lie in $P\cup\{a,b\}$.  Otherwise $P$ is common to all members.
For (ii), take incident edges $ab,ac$.  An edge avoiding $a$ must be
$bc$ in order to meet both; any other edge meeting all three of
$ab,ac,bc$ is also one of those three edges.  Cases with at most one
edge are immediate.
\end{proof}

The Erd\H{o}s--Gallai matching theorem \cite{EG59} will be used in the
form: a simple graph on $N\ge7$ vertices with no matching of size three
has at most $2N-3$ edges.

\begin{lemma}[No common two-set]\label{lem:nocore}
Assume $n\ge14$.  If $\mathcal H$ is $2$-intersecting, has no common
two-set, and contains two supports meeting in exactly two coordinates,
then
\begin{equation}\label{eq:nocore}
 |C_4|\le4n+15.
\end{equation}
\end{lemma}
\begin{proof}
Choose $S,T\in\mathcal H$ with $P=S\cap T$ of size two and write
$S=P\sqcup A$, $T=P\sqcup B$, where $P,A,B$ are pairwise disjoint
two-sets.  Since $P$ is not common, choose $R\in\mathcal H$ not
containing $P$.  Split the supports into
\[
 \mathcal H_0=\{U\in\mathcal H:P\not\subseteq U\},\qquad
 \mathcal H_1=\{U\in\mathcal H:P\subseteq U\}.
\]
For any subfamily $\mathcal F\subseteq\mathcal H$, denote its
\emph{excess} over one word per support by
$\sum_{U\in\mathcal F}(\mu(U)-1)$.

\emph{First, the excess on $\mathcal H_0$ is at most eight.}
If $U\in\mathcal H_0$ has multiplicity greater than one, then
$2\le |U\cap S|,|U\cap T|\le3$ by the two-intersection
property and the fact that $U$ is distinct from $S,T$.  Moreover,
\[
 |U\cap S|+|U\cap T|
 =|U\cap(S\cup T)|+|U\cap P|\le5.
\]
The two intersection sizes cannot both be two: in that case the
corresponding exact-overlap pairs are different, and
Lemma~\ref{lem:overlaps} would force $\mu(U)=1$.
Consequently their sizes are two and three in some order.
The supports of type $(2,3)$ are precisely
\[
 U=\{p,a\}\cup B\quad (p\in P,\ a\in A),
\]
and the supports of type $(3,2)$ are
$U=\{p,b\}\cup A$ with $p\in P,b\in B$.
Within either group of four supports, at most two can have multiplicity
greater than one.  Indeed, any three pairs $(p,a)$ chosen from a
$2\times2$ grid include two with different first and second indices.
Their supports intersect in exactly $B$, respectively $A$, and this
gives a second exact-overlap pair, forcing multiplicity one.
Each of the at most four remaining supports has multiplicity at most
three by Lemma~\ref{lem:overlaps}.  Hence
\begin{equation}\label{eq:H0excess}
 \sum_{U\in\mathcal H_0}(\mu(U)-1)\le8.
\end{equation}

\emph{Next, bound the excess on $\mathcal H_1$.}
Put $r=|R\cap P|\in\{0,1\}$.
If $r=0$, every $U\in\mathcal H_1$ equals $P\cup e$ with
a two-set $e\subseteq R$.  There are at most six such supports,
and each meets $R$ in exactly two positions.  Thus each has
multiplicity at most three and the excess on $\mathcal H_1$ is
at most $12$.

Suppose $r=1$.  Write $R\cap P=\{p\}$ and $D=R\setminus P$;
then $|D|=3$.  Represent $\mathcal H_1$ by a graph $G$ on the
coordinates outside $P$, with an edge $e$ for each support $P\cup e$.
Every edge of $G$ meets $D$.  Call an edge \emph{exceptional} if it
is contained in $D$; there are at most three.  Every other edge is
\emph{ordinary} and its support meets $R$ in exactly two coordinates.
An ordinary edge therefore has multiplicity at most three.

If $G$ has two disjoint edges, an ordinary edge with multiplicity
greater than one must meet every edge of $G$.  Otherwise it would
have two different exact-overlap pairs, one from $R$ and one from its
disjoint external edge, contradicting Lemma~\ref{lem:overlaps}.
Relative to two fixed disjoint edges, any edge meeting both is a
cross-edge in a $K_{2,2}$.  At most two of those cross-edges can
meet every edge, because the universal cross-edges are themselves
pairwise intersecting.  Their total excess is therefore at most
$2(3-1)=4$.
If $k\le3$ exceptional edges occur, the corresponding supports
intersect pairwise in three positions.  Their total multiplicity is
at most $9$ for $k=1$, at most $15$ for $k=2$ by
Lemma~\ref{lem:common3}, and at most $21$ for $k=3$ by
Corollary~\ref{cor:three-small}.  Their excess is therefore at most
$18$ in every case.  Thus the excess on $\mathcal H_1$ is at most
$4+18=22$.

If $G$ has no two disjoint edges, Lemma~\ref{lem:structure}(ii)
says that $G$ is contained in a triangle or is a star.
For a triangle, at most three supports occur in $\mathcal H_1$;
they pairwise intersect in three positions, so their total
multiplicity is at most $21$.  Therefore their excess is at most
$21$.  If $G$ is a star centered outside $D$, all its edges must
lead into $D$, so it has at most three ordinary edges; their excess
is at most $3(3-1)=6$.

Finally, suppose $G$ is a star centered at $c\in D$.
At most two star edges are exceptional.  By
Lemma~\ref{lem:common3}, their total multiplicity is at most $15$.
For every ordinary edge $\{c,v\}$ with $v\notin D$, comparison
with a word on $R$ fixes the two symbols on the same pair $\{p,c\}$.
If two different ordinary leaves $v,w$ exist, comparison of their
words forces the symbol at the other coordinate of $P$ to agree:
otherwise the two leaf coordinates contribute two and that
additional disagreement gives the forbidden distance three.
Thus all three common-core symbols agree across ordinary leaves,
and each ordinary support has multiplicity one; two words on it
could otherwise differ only at the leaf coordinate.
If at most one ordinary leaf exists, its multiplicity is at most
three.  Hence if $t$ ordinary edges occur, their total multiplicity
is at most $t+2$.  Adding at most $15$ words on exceptional edges
shows that the excess on $\mathcal H_1$ is at most $17$.

In every case, the excess on $\mathcal H_1$ is at most $22$.
Together with \eqref{eq:H0excess}, the complete nontrivial
intersection bound \eqref{eq:AK} now yields
\[
 |C_4|\le|\mathcal H|+8+22
 \le(4n-15)+30=4n+15.
\]
\end{proof}

\section{When a common two-set exists}

Assume all supports in $\mathcal H$ contain a fixed two-set $P$.
Associate a graph $G$ on $[n]\setminus P$ by putting
$e\in E(G)$ exactly when $P\cup e\in\mathcal H$.
The \emph{matching number} of $G$ is the largest number of pairwise
disjoint edges in $G$.

\begin{lemma}[Three disjoint external edges]\label{lem:matching3}
If $G$ has three pairwise disjoint edges and $n\ge6$, then
\begin{equation}\label{eq:matching3}
 |C_4|\le\binom{n-2}{2},\qquad |C_2|\le2n-3.
\end{equation}
\end{lemma}
\begin{proof}
Let $e_1,e_2,e_3$ be disjoint edges of $G$ and choose words on
$P\cup e_i$.  Any two of these supports intersect only in $P$,
so their external coordinates already contribute distance four.
Consequently the three words agree on $P$.  An arbitrary word in
$C_4$ has an external pair disjoint from at least one of the three
edges, hence comparison with that word forces the same symbols on $P$.
Deleting $P$ is thus an isometry from $C_4$ to a pure weight-two code
of length $n-2$.  Equation~\eqref{eq:pure-w2} gives the first bound.

Every weight-two support must meet each of $P\cup e_i$, since a
disjoint two-support would produce distance six.  A two-set avoiding
$P$ cannot meet all three disjoint $e_i$, so every weight-two
support meets $P$.  Label $P=\{a,b\}$.  The possibilities are
$P$, $\{a,v\}$ and $\{b,v\}$ for $v\notin P$.
For $\{a,v\}$, compare with a four-word whose external pair avoids
$v$: Lemma~\ref{lem:one} fixes the symbol at $a$.  Thus this
support has multiplicity at most one, and similarly for $\{b,v\}$.
If any such noncore support occurs, comparing with it fixes one
coordinate of every weight-two word on $P$, so that $P$ also has
multiplicity at most one.  This gives $2(n-2)+1=2n-3$ words.
If no noncore support occurs, the multiplicity on $P$ is at most
three, also at most $2n-3$.  This proves the second bound.
\end{proof}

\begin{lemma}[A matching of size two]\label{lem:matching2}
Suppose $G$ contains two disjoint edges but no matching of size three.
For $n\ge14$,
\begin{equation}\label{eq:matching2}
 |C_4|\le2n+12.
\end{equation}
\end{lemma}
\begin{proof}
Choose disjoint edges $e,f$.  Let $X$ be the set of edges of $G$
meeting both $e$ and $f$.  It is contained in the $K_{2,2}$ between
the endpoint sets of $e$ and $f$, so $|X|\le4$.
All codewords corresponding to edges outside $X$ have the same two
symbols on $P$: the words on $e$ and $f$ agree on $P$, and each
other edge outside $X$ is disjoint from $e$ or $f$.
Delete $P$ from this part of $C_4$.  We obtain a pure weight-two
code whose distinct supports are the edges of $E(G)\setminus X$.
Its support graph has at most two isolated edges, since three isolated
edges would form a matching of size three in $G$.
The first inequality in \eqref{eq:pure-w2} bounds this part by
$|E(G)\setminus X|+4$ words.

The supports corresponding to $X$ contribute at most $15$ words.
Indeed, with zero or one edge this follows from
Lemma~\ref{lem:local9}.  Two opposite edges have disjoint external
parts, so each has multiplicity at most three; two adjacent edges
have a common three-set, so Lemma~\ref{lem:common3} bounds their
combined multiplicity by $15$.  If three cross-edges occur, the
two opposite edges have multiplicity at most three and the third
at most nine, giving $15$.  If all four cross-edges occur, each
has a disjoint opposite edge and hence multiplicity at most three,
giving at most $12$.

The Erd\H{o}s--Gallai theorem gives
$|E(G)|\le2(n-2)-3=2n-7$.  Consequently
\[
 |C_4|\le |E(G)\setminus X|+4+15
 \le (2n-7)+19=2n+12.
\]
\end{proof}

\begin{lemma}[No two disjoint external edges]\label{lem:matching1}
Suppose $G$ has no two disjoint edges.  If $G$ lies in a triangle,
then $|C_4|\le27$.  If $G$ is a star with at most two edges, then
$|C_4|\le18$.  If $G$ is a star with at least three edges and
$n\ge14$, then
\begin{equation}\label{eq:matching1}
 |C_4|\le4n-12,\qquad |C_2|\le3n.
\end{equation}
\end{lemma}
\begin{proof}
The first two assertions follow from
Lemma~\ref{lem:local9}: there are at most three or two supports,
respectively.  In the remaining case every external edge contains
a common vertex $c$, so all four-supports contain the same
three-set $Q=P\cup\{c\}$.  If there are $s$ distinct supports,
then $3\le s\le n-3$.  By Lemma~\ref{lem:common3},
\[
 |C_4|\le\max\{4s,2s+11\}\le4n-12
\]
for $n\ge14$.

There are at least three distinct leaf coordinates.  A two-set
avoiding $Q$ meets at most two of the four-supports and would
therefore be disjoint from a third, which is impossible for a
weight-two support.  Thus every weight-two support meets $Q$.
There are $3(n-3)$ supports meeting $Q$ in exactly one coordinate.
On any such support $\{i,v\}$, comparison with a four-word on a
leaf other than $v$ fixes the symbol at $i$ and hence forces
multiplicity one.  The three two-supports lying inside $Q$ have
multiplicity at most three each.  Therefore
$|C_2|\le3(n-3)+9=3n$.
\end{proof}

\section{The combinatorial range \texorpdfstring{$n\ge16$}{n >= 16}}

\begin{theorem}\label{thm:large}
For every $n\ge16$, a quaternary $\{2,4\}$-code has size at most
$1+\binom n2$.
\end{theorem}
\begin{proof}
Normalize as in \eqref{eq:layers}.  If $C_4$ is empty, use
\eqref{eq:pure-w2}.  Otherwise \eqref{eq:C2-general} holds,
and the support family $\mathcal H$ is $2$-intersecting by
Lemma~\ref{lem:two-intersecting}.

Suppose first that $\mathcal H$ has no common two-set.
If every two supports intersect in at least three coordinates,
Lemma~\ref{lem:structure}(i) and the absence of a common three-set
show that all supports lie in a five-set.  Corollary~\ref{cor:three-small}
and \eqref{eq:C2-general} give
\begin{equation}\label{eq:case1a}
 |C|\le1+(4n-2)+33=4n+32.
\end{equation}
Otherwise two supports intersect in precisely two coordinates.
Lemmas~\ref{lem:C2-bounds} and \ref{lem:nocore} give
\begin{equation}\label{eq:case1b}
 |C|\le1+(2n+9)+(4n+15)=6n+25.
\end{equation}
For $n\ge16$, both \eqref{eq:case1a} and \eqref{eq:case1b}
are at most $1+\binom n2$.  Indeed, for the second and closer
estimate,
\[
 \binom n2-(6n+24)=\frac{(n-16)(n+3)}2\ge0;
\]
the first estimate is already strictly smaller at $n=16$ and has
smaller linear growth.

Now suppose every member of $\mathcal H$ contains a common two-set,
and form its external graph $G$.  If $G$ has a matching of size
three, Lemma~\ref{lem:matching3} gives
\[
 |C|\le1+(2n-3)+\binom{n-2}{2}=1+\binom n2.
\]
If $G$ has a matching of size two but none of size three, then
two four-supports meet in exactly two coordinates, and
Lemmas~\ref{lem:C2-bounds} and \ref{lem:matching2} give
\[
 |C|\le1+(2n+9)+(2n+12)=4n+22.
\]
Finally, if $G$ has no two disjoint edges, use
Lemma~\ref{lem:matching1}.  In the triangle case,
$|C|\le1+(4n-2)+27=4n+26$; for a star with at most two edges,
$|C|\le4n+17$; and for a star with at least three edges,
$|C|\le1+3n+(4n-12)=7n-11$.
All of these are at most $1+\binom n2$ at $n=16$, and their
linear increments are strictly smaller than the successive
increments $n$ of $\binom n2$ for $n\ge16$.
Every structural possibility has been covered.
\end{proof}

\section{Finite verification for \texorpdfstring{$6\le n\le15$}{6 <= n <= 15}}
\label{sec:finite}

The following is an exact finite proof, not a numerical optimization
estimate.  Two supplementary C++17 programs implement the
mathematical procedures below using integer representations and
exhaustive coloring branch-and-bound.  All reported runs were
executed for this version of the manuscript.

\subsection{The normalized candidate graph}

Suppose $C$ were a counterexample.  By \eqref{eq:pure-w2},
its weight-four layer is nonempty.  Fix one four-word and apply a
Hamming isometry fixing zero to normalize it to
\[
 z=(1,1,1,1,0,\ldots,0).
\]
Let $V_n$ consist of all $x\ne z$ with
$\wt(x)\in\{2,4\}$ and $d(x,z)\in\{2,4\}$.
Make a graph $\Gamma_n$ on $V_n$ by joining distinct words whose
Hamming distance is two or four.  The code minus $\{0,z\}$ is then
a clique in $\Gamma_n$.  Thus a counterexample to
$|C|\le1+\binom n2$ would give a clique of size
at least $\binom n2$ in $\Gamma_n$.

The candidate count can be computed exactly.  For weight two,
the support meets the first four coordinates in one or two places,
yielding respectively $12(n-4)$ and $30$ choices.
For weight four, the intersection size is two, three, or four,
yielding respectively $54\binom{n-4}{2}$,
$156(n-4)$, and $40$ choices.  Hence
\begin{equation}\label{eq:candidate-count}
 |V_n|=12n-18+54\binom{n-4}{2}+156(n-4)+40
 =27n^2-75n-62.
\end{equation}
To see the intermediate factors: for support overlap two, both
common symbols must be $1$; the two other coordinates each have
three nonzero choices.  For overlap three, there are
$1+\binom32\cdot2^2=13$ allowed common-symbol patterns,
three choices for the outside nonzero symbol, and four choices of
common support.  For overlap four, exactly two or four of the four
symbols differ from $1$, giving $\binom42 2^2+2^4=40$.
The weight-two count follows by the same distance formula.

\subsection{The range \texorpdfstring{$6\le n\le13$}{6 <= n <= 13}}

The group fixing both $0$ and $z$ permutes the first four
coordinates, permutes the remaining ones, fixes $0$ and $1$
on the first four coordinates while possibly exchanging $2,3$,
and permutes the three nonzero symbols independently in the
other coordinates.  Its orbits on $V_n$ are classified by the triple
\[
 \bigl(\wt(x),\ |\supp(x)\cap\supp(z)|,\
 |\{i\in\supp(x)\cap\supp(z):x_i\ne1\}|\bigr).
\]
There are exactly eight orbits, with triples
\begin{equation}\label{eq:eight-orbits}
 (2,1,0),\ (2,2,0),\ (2,2,2),\
 (4,2,0),\ (4,3,0),\ (4,3,2),\ (4,4,2),\ (4,4,4).
\end{equation}
Fix one representative $v$ from each orbit.  If $\Gamma_n$ contains
a clique of size $\binom n2$, then for at least one representative
$v$ the induced graph on its neighbors contains a clique of size
$\binom n2-1$.  This is an equivalence, since the stabilizer maps
any selected clique vertex to the representative of its orbit.

Our exact verifier constructs $V_n$ from its definition, verifies
\eqref{eq:candidate-count}, enumerates the eight triples
\eqref{eq:eight-orbits}, and tests all eight forbidden clique sizes.
It returns \texttt{NO} in every case.  The total recursion counts
in Table~\ref{tab:small} are reproducibility data, not assumptions
needed for the correctness proof.

\begin{table}[ht]
\centering
\caption{Independent exact clique decisions for the eight orbit
representatives.  In every row the forbidden clique size in
$\Gamma_n$ is $\binom n2$.}
\label{tab:small}
\begin{tabular}{rrrr}
\toprule
$n$ & $|V_n|$ & forbidden size & total recursive nodes \\
\midrule
6 & 460 & 15 & 2141 \\
7 & 736 & 21 & 3416 \\
8 & 1066 & 28 & 1446 \\
9 & 1450 & 36 & 2055 \\
10 & 1888 & 45 & 2888 \\
11 & 2380 & 55 & 6752 \\
12 & 2926 & 66 & 29437 \\
13 & 3526 & 78 & 130166 \\
\bottomrule
\end{tabular}
\end{table}

\subsection{The range \texorpdfstring{$n=14,15$}{n = 14, 15}}

The structural estimates already cover every case except when
$\mathcal H$ has no common two-set and two supports meet in exactly
two coordinates.  Indeed, if every two supports meet in three
coordinates, \eqref{eq:case1a} applies; a common-two-set family is
covered by Lemmas~\ref{lem:matching3}--\ref{lem:matching1}.
Substituting $n=14,15$ into those estimates gives the desired
upper bound in all such cases.

For the remaining case, relabel coordinates so that two
four-supports are
\[
 S=\{0,1,2,3\},\qquad T=\{0,1,4,5\},\qquad
P=S\cap T=\{0,1\}.
\]
(Here coordinates are labeled $0,\ldots,n-1$ only for this finite
calculation.)  Every occurring four-support $U$ meets both $S$
and $T$ in at least two positions.  The number of possible
four-supports satisfying these two conditions is
\begin{equation}\label{eq:support-count}
 \binom{n-2}{2}+8(n-6)+9=\frac{n^2+11n-72}{2},
\end{equation}
namely $139$ for $n=14$ and $159$ for $n=15$.
To obtain this count, a candidate containing $P$ is determined by
any two of the remaining $n-2$ coordinates.  If it contains exactly
one point of $P$, there are $8(n-6)+8$ possibilities; if it avoids
$P$, it must be $(S\setminus P)\cup(T\setminus P)$,
which gives the last one.

For each candidate four-support $U$, introduce eight abstract
\emph{states}:
\begin{enumerate}[label=\textup{(\alph*)}]
\item a \emph{light} state of weight $1$;
\item six \emph{heavy} states, one for each two-subset $L\subset U$,
 each of weight $3$;
\item a \emph{free} state of weight $9$.
\end{enumerate}
States belonging to the same support are not adjacent.  States
on distinct supports $U,V$ are adjacent precisely when
$|U\cap V|\ge2$ and, if $|U\cap V|=2$, neither state is free and
any heavy state must be indexed by the pair $U\cap V$.
This defines a finite weighted graph.

\begin{lemma}[Soundness of the support-state relaxation]
\label{lem:state}
For every actual code with the fixed supports $S,T$, the supports
of its weight-four layer admit states forming a clique in the
weighted graph, with total state weight at least $|C_4|$.
The state of $S$ and the state of $T$ can each be restricted to the
light state or the heavy state indexed by $P$.
\end{lemma}
\begin{proof}
Consider any occurring support $U$.  If it has no exact two-overlap
with another occurring support, assign its free state;
Lemma~\ref{lem:local9} bounds its multiplicity by nine.
If all exact-overlap pairs $U\cap V$ are the same two-set $L$,
assign the heavy state indexed by $L$;
Lemma~\ref{lem:overlaps} bounds its multiplicity by three.
If two different exact-overlap pairs occur, assign its light
state; the second assertion of Lemma~\ref{lem:overlaps} forces
multiplicity one.
All pairs of assigned states are adjacent by construction,
because exact intersections exclude a free state and force any
heavy index to equal the intersection.  The assigned state weights
dominate the corresponding multiplicities, proving the first claim.
Since $S$ and $T$ have the exact-overlap pair $P$, neither can be
free, and any heavy state on either support must use $P$.
\end{proof}

Since $P$ is not common, some third support $R$ avoids containing
$P$.  Up to permutations fixing the ordered pair $(S,T)$, the
possible intersection types of $R$ have representatives
\[
 \{2,3,4,5\},\quad \{0,2,4,6\},\quad
 \{0,2,4,5\},\quad \{0,2,3,4\}.
\]
Indeed, write
\[
 (p,a,b,o)=\bigl(|R\cap P|,|R\cap(S\setminus P)|,
 |R\cap(T\setminus P)|,|R\setminus(S\cup T)|\bigr).
\]
The conditions $p\le1$, $p+a\ge2$, $p+b\ge2$ and
$p+a+b+o=4$ give exactly
$(0,2,2,0)$, $(1,1,1,1)$, $(1,1,2,0)$ and $(1,2,1,0)$.
After the two allowed choices for each of the states of $S,T$,
and compatibility of the third state, there are respectively
$1,1,4,4$ forced-state choices in these four types: ten cases
in total.

A complete weighted-clique decision procedure on these ten cases
excludes a total state weight of at least $55$ for $n=14$ and at
least $67$ for $n=15$.  These are exactly the thresholds required
for the theorem: Lemma~\ref{lem:C2-bounds} gives
$|C_2|\le37$ at $n=14$ and $|C_2|\le39$ at $n=15$.
A counterexample of cardinality at least $93$, respectively $107$,
would therefore require $|C_4|\ge55$, respectively $67$.
The reproduced integer-only search statistics are
\begin{center}
\begin{tabular}{rrrrr}
\toprule
$n$ & supports & abstract states & forced cases & recursive nodes\\
\midrule
14 & 139 & 1112 & 10 & 58\\
15 & 159 & 1272 & 10 & 27\\
\bottomrule
\end{tabular}
\end{center}
The entries report a \texttt{NO} decision for every forced case.

\subsection{Why the exact searches are exhaustive}
\label{sec:algorithm}

For an ordinary or weighted graph, let a recursive node contain a
clique already selected and a set of vertices compatible with all
selected vertices.  Properly color the induced candidate graph by
a deterministic greedy procedure, obtaining independent classes
$D_1,\ldots,D_k$.  For ordinary cliques, a clique extension uses at
most $k$ vertices.  For weighted cliques, an extension has weight
at most
\begin{equation}\label{eq:coloring-bound}
 \sum_{i=1}^k\max_{v\in D_i}w(v).
\end{equation}
The algorithm branches on a chosen vertex $v$ by including $v$ and
intersecting the remaining candidates with its neighborhood, or by
excluding $v$.  These two alternatives partition all possible
extensions.  A branch is discarded only when the weight already
selected plus the proper-coloring upper bound is strictly below the
forbidden threshold.  Induction on the number of remaining candidates
proves that this returns \texttt{NO} if and only if the forbidden
clique cannot occur.

Both supplied implementations strengthen this bound by ordering
vertices in successive color classes and using the valid bound
for every prefix of that order.  This is a standard exact
branch-and-bound optimization, not a change in the decision problem.
All vertex generation, distance comparisons, colorings, weights,
and branch decisions use integers.  In particular, the state graph
is only a \emph{relaxation}: its negative decision proves the required
bound for real codeword configurations without presupposing that every
state clique is realizable.

\begin{theorem}[Finite range]\label{thm:finite}
For $6\le n\le15$, every quaternary $\{2,4\}$-code has size at most
$1+\binom n2$.
\end{theorem}
\begin{proof}
The pure weight-two case follows from \eqref{eq:pure-w2}.
For $6\le n\le13$, the normalized candidate graph and all eight
orbits are covered by the exact decisions of Table~\ref{tab:small}.
For $n=14,15$, the structural bounds cover all but the no-common-core
exact-overlap case, and Lemma~\ref{lem:state} reduces that case to the
ten weighted decisions just described.  Their negative answers
exclude any counterexample.  In all cases the decisions are exhaustive
by the preceding correctness argument.
\end{proof}

\section{Completion and reproducibility}

\begin{proof}[Proof of Theorem~\ref{thm:main}]
The standard construction \eqref{eq:standard} gives the lower
bound.  Theorems~\ref{thm:large} and \ref{thm:finite} give the
quaternary upper bound for every $n\ge6$.  Finally,
\eqref{eq:alphabet-monotone} gives the ternary upper bound,
and the same construction attains equality for $q=3$.
\end{proof}

\noindent\textbf{Data and code availability.}
The two supplementary programs are
\begin{center}
\texttt{q34\_small\_exact.cpp}\qquad
\texttt{q34\_weighted\_exact.cpp}
\end{center}
They implement the independent finite decision reductions using
standard C++17 and no external libraries.  Exact console outputs
for the eight small dimensions and two weighted dimensions are
included as supplementary logs.  A reproducible compiler command is
\begin{center}
\texttt{g++ -O3 -std=c++17 SOURCE.cpp -o verifier}
\end{center}
Run the first executable with each integer argument $6,\ldots,13$
and the second with $14$ or $15$.
The source files regenerate all graph vertices and adjacencies from
their mathematical definitions; there are no hard-coded clique
certificates or external data files.

\medskip
\noindent\textbf{Declaration of AI-assisted preparation.}
OpenAI ChatGPT was used to assist with manuscript organization,
editing, mathematical cross-checking, and construction of
independent verification software.  The authors are responsible for
all mathematical claims, program outputs, citations, and final
submission.  The exact verifiers do not invoke an AI system.

\end{document}